\documentclass[12pt]{article}
\usepackage[all]{xy}

\usepackage{epic}
\usepackage{amsmath}[1996/11/01]
\usepackage{amssymb,amsthm,amsfonts,latexsym,epsfig,graphics}
 \def\beql#1#2\eeql{\begin{equation}\label{#1}#2\end{equation}}

\advance \textheight by -1 \topmargin
\advance \textheight by -1 \headheight
\advance \textheight by -1 \headsep
\advance \textheight by -1 \footskip
\vsize = \textheight
\hsize = \textwidth
\newcommand{\operp}{\mathop{\bigcirc\!\!\!\!\!\!\!\perp\,}} %orthog. dir. Summe
\DeclareMathOperator{\Syl}{Syl}

\DeclareMathOperator{\Stab}{Stab}

\DeclareMathOperator{\GL}{GL}
\DeclareMathOperator{\Aut}{Aut}
\DeclareMathOperator{\Alt}{Alt}

\DeclareMathOperator{\im}{im}

\DeclareMathOperator{\SL}{SL}
\DeclareMathOperator{\PSL}{PSL}

\newtheorem{theorem}{Theorem}[section]

\newcommand{\bew}{\noindent\underline{Proof.}\ }

\newtheorem{remark}[theorem]{Remark}

\newtheorem{proposition}[theorem]{Proposition}

\newtheorem{defn}[theorem]{Definition}

\newcommand{\Z}{{\mathbb{Z}}}
\newcommand{\Q}{{\mathbb{Q}}}

\newcommand{\R}{{\mathbb{R}}}

\newcommand{\trace}{\mbox{trace}}

\newcommand{\eb}{\phantom{zzz}\hfill{$\square $}\smallskip}

\begin{document}
\begin{center}
{\Large {\bf A fourth
 extremal even unimodular lattice of dimension 48.}}\\ 
\vspace{1.5\baselineskip}
{\em Gabriele Nebe} \\
\vspace*{1\baselineskip}
Lehrstuhl D f\"ur Mathematik, RWTH Aachen University\\
52056 Aachen, Germany \\
 nebe@math.rwth-aachen.de \\
\vspace{1.5\baselineskip}
\end{center}

{\sc Abstract.}
{\small 
We show that there is a unique 
extremal even unimodular lattice of dimension 48 which has an automorphism of
order 5 of type $5-(8,16)-8$. Since the three known extremal lattices 
do not admit such an automorphism, this provides a new example of
an extremal even unimodular lattice in dimension 48.
}

\section{Introduction}

A lattice $L$ in Euclidean space $(\R^n,(,))$ is a
 free $\Z $-module of
 rank $n$ containing a basis of $\R^n$. 
The lattice $L$ is called {\em even}, if the 
associated quadratic form $Q:\R^n \to \R_{\geq 0} , x\mapsto Q(x) := \frac{1}{2} (x,x) $ is integral on $L$, so $Q(L) \subseteq \Z $. 
Then $L$ is contained in its dual lattice 
$$L^{\#} := \{ y\in \R^n \mid (y,\ell )\in \Z \mbox{ for all } \ell \in L \} .$$ 
$L$ is {\em unimodular}, if $L=L^{\#}$. 
Any lattice $L$ defines a sphere packing,
whose density measures its error correcting properties.
One of the main goals in lattice theory is to find dense lattices.
This is a very difficult problem, the densest lattices are known only in dimension $n\leq 8$
and in dimension $24$ \cite{CohnKumar}, for $n=8$ and $n=24$ the densest lattices
are even unimodular lattices.
The density of a unimodular lattice is measured by its {\em minimum},
$$\min (L):= \min \{ 2 Q(\ell ) \mid  0\neq \ell \in L \}.$$
For even unimodular lattices the theory of modular forms allows one to bound
this minimum $\min (L) \leq 2+ 2 \lfloor \frac{n}{24} \rfloor$
and {\em extremal lattices} are those even unimodular lattices $L$ that
achieve equality.
Of particular interest are
extremal even unimodular lattices $L$ in the {\em jump dimensions} $24m$.
For $m=1$ there is a  unique extremal even unimodular lattice,
the {\em Leech lattice}, which is the densest
 24-dimensional lattice \cite{CohnKumar}.
By \cite{Conway}, its
automorphism group 
%$$\Aut(\Lambda _{24})  =\{ g\in \GL (\Lambda _{24}) \mid 
%Q(xg) = Q(x) \mbox{ for all } x \}  \cong 2.Co_1 $$
is a covering group of the
sporadic simple Conway group $Co_1$. 
The 196560 minimal vectors of the Leech lattice form the
unique tight spherical 11-design and realise the maximal
kissing number in dimension 24.
In dimension 72 one knows one extremal unimodular lattice
\cite{DIM72}. The existence of such a lattice was a longstanding 
open problem. 
In dimension 48 there are at least four extremal 
even unimodular lattices. They 
are the densest known lattices in their dimension and
realise the maximal known kissing number $52,416,000$.
It is a very interesting problem to classify all 
48-dimensional extremal even unimodular lattices. 
To get an idea of how many such lattices might exist, I
started a program to find all extremal lattices $L$ 
whose automorphism group
$$
\Aut(L)  =\{ g\in \GL (L) \mid 
Q(xg) = Q(x) \mbox{ for all } x \} $$
is not too small. 
In  \cite{DIM48} I classified all 48-dimensional extremal
lattices that have an automorphism of order $a$ whose Euler phi value
is $\varphi(a) > 24$. All these lattices 
are isometric to one of the lattices $P_{48p}$,
$P_{48q}$, or $P_{48n}$, which were known before. 
The present paper classifies all extremal lattices invariant 
under a certain automorphism of order 5. 
It turns out that there is a unique such lattice, $P_{48m}$, and this 
lattice is not isometric to one of the lattices above. 

\begin{table}[ht]\label{t2} 
{
\begin{tabular}{|lccr|}
\hline
\multicolumn{4}{|c|}{Table \ref{t2}: The known extremal even unimodular lattices in the jump dimensions} \\
\hline
\hline
name & autom. group & order & ref. \\ \hline
$\Lambda _{24} $ &  $2.Co_1$ &  $ 8315553613086720000 $ &  \\
		 &   & $ =  2^{22} 3^9 5^47^2 11\ 13\ 23 $  & \cite{Leech}, \cite{Conway} \\
\hline
$P_{48p}  $ &  $(\SL _2(23) \times S_3) : 2$ &   $72864 = 2^53^211\ 23  $ 
& \cite{SPLAG}, \cite{DIM48} \\
\hline
 $P_{48q} $  & $  \SL _2(47) $ &   $ 103776 =2^53\ 23\ 47 $  & \cite{SPLAG}, \cite{DIM48}\\
\hline
$P_{48n}  $ &  $ (\SL _2(13) {\sf Y} \SL_2(5) ).2^2 $ &   $ 524160 = 2^73^25\ 7\ 13 $ & \cite{cyclo}, \cite{DIM48}  \\
\hline
$P_{48m} $  & & multiple of  & this \\
& &    $ 1200 = 2^43\ 5^2 $ &  paper \\
\hline
$\Gamma _{72}  $ &  $(\SL _2(25) \times \PSL _2(7) ) : 2$ &  $5241600= 2^8  3^2  5^2 7\ 13  $ & \cite{DIM72}, \cite{DIM48} \\
\hline
\end{tabular}}
\end{table}

\section{The type of an automorphism}

The notion of the type of an automorphism of a lattice $L$ 
was introduced in \cite{DIM48}. It was motivated by the 
analogous notion of a type of an automorphism of a code. 

Let $\sigma \in \GL_n(\Q )$ be an element of prime order $p$. 
Let $K:=\ker (\sigma -1)$  and $I:=\im (\sigma -1)$.
Then $K$ is the fixed space of $\sigma $ and the action of 
$\sigma  $ on $I$ gives rise to a vector space structure on $I$ over the 
$p$-th cyclotomic number field $\Q[\zeta _p]$.
In particular $n= d+z(p-1)$, where 
$d:=\dim _{\Q} (K)$  and $z = \dim _{\Q[\zeta _p]} (I)$. 

If $L$ is a $\sigma $-invariant $\Z$-lattice, then 
$L$ contains a sublattice  $M$ with
$$ L \geq M= (L\cap K) \oplus (L \cap I) = : L_K(\sigma ) \oplus L_I (\sigma ) \geq pL $$ 
of finite index $[L:M] = p^s$ where $s\leq \min (d,z )$.

\begin{defn}
The triple $p-(z,d)-s$ is called the {\em type} of the element $\sigma \in \GL(L)$.
\end{defn}

\begin{remark}\label{unimod}
Let $(L,Q)$ be an even unimodular lattice and $\sigma \in \Aut (L) $
be of type $p-(z,d)-s$. Then 
$L_I^{\#}(\sigma )/L_I(\sigma ) \cong (\Z[\zeta _p]/(1-\zeta _p)) ^s$ and 
$L_K^{\#}(\sigma )/L_K(\sigma ) \cong (\Z/p\Z)^s$ as $\Z[\sigma ]$-modules.
In particular $0\leq s\leq \min (z,d )$.
If $z=s$ then $(1-\sigma) L_I^{\#} (\sigma ) = L_I(\sigma) $
and hence $(L_I^{\#}(\sigma ), pQ ) = ( X, \trace_{\Q[\zeta _p]/\Q }(h))$ 
is the trace lattice 
of an Hermitian unimodular $\Z[\zeta _p]$ lattice $(X,h)$ of rank 
$\dim _{\Z[\zeta_p]} (X) = z$. 
\end{remark}

\begin{remark} \label{line12}
If $L$ is an even lattice and $p$ is odd, then $L_K(\sigma )$ and $L_I(\sigma )$ are also
even lattices, because $L_K(\sigma ) \oplus L_I (\sigma )$ is a sublattice of 
odd index in $L$. 
\end{remark} 

In \cite{DIM48} we narrowed down the possible types of 
prime order automorphisms of an extremal even unimodular lattice in dimension 48. 
By Remark \ref{line12} the fixed lattice of an element of order 3
 cannot be $\sqrt{3} D_{12}^+$, as this is an odd lattice. So Type 
$3-(18,12)-12$ is not possible and the possible types
are among the ones in Table \ref{t1}. 

\begin{table}[ht]\label{t1} 
{\begin{tabular}{|c|c|c|c|l|}
\hline
\multicolumn{5}{|c|}{Table \ref{t1}: 
The possible types of automorphisms $\sigma \neq -1$ 
of prime order} \\
\hline
\hline
 Type & $L_K (\sigma ) $ &  $L_I  (\sigma ) $ & example & complete\\
\hline
 47-(1,2)-1 &  unique &  unique & $P_{48q}$ & \cite[Thm 5.6]{DIM48} \\
\hline
 23-(2,4)-2  & unique & at least 2 & $P_{48q}$, $P_{48p}$ &  \\
\hline
 13-(4,0)-0 &  $\{ 0 \}$ & at least 1  &  $P_{48n}$ &  \\
\hline
 11-(4,8)-4 &  unique & at least 1 & $P_{48p}$ &  \\
\hline
 7-(8,0)-0 &  $\{ 0 \}$ & at least 1 & $P_{48n}$ &   \\
 7-(7,6)-5 & $\sqrt{7}A_6^{\#}$ & not known & not known &  \\
\hline
 5-(12,0)-0 &  $\{0 \} $ & at least 2 & $P_{48n}$, $P_{48m}$ &   \\
 5-(10,8)-8 & $\sqrt{5}E_8$ & at least 1 &   $P_{48m}$  &  \\
 5-(8,16)-8 & $[2.\Alt_{10}]_{16}$ & $\Lambda _{32}$ & $P_{48m}$ & Thm. \ref{line9}  \\
\hline
 3-(24,0)-0 & $\{0 \}$  & at least 3 & $P_{48p}$, $P_{48n}$,  $P_{48m}$ &  \\
 3-(20,8)-8 & $\sqrt{3} E_8$ & not known & not known &  \\
 3-(16,16)-16 &  $\sqrt{3} (E_8\perp E_8) $ & at least 4 & 
 $P_{48p}$, $P_{48q}$, $P_{48n}$ &   \\
 3-(16,16)-16 & $\sqrt{3} D_{16}^+$ & at least 4 & not known &  \\ 
 3-(15,18)-15 &  unique & two & not known &   \\
3-(14,20)-14 &  ?  &  unique & not known &   \\
 3-(13,22)-13 & ? &  unique & not known &  \\
\hline
 % 2-(48,0)-0 & $\{ 0 \}$ & \multicolumn{3}{|c|}{$\sigma = -1$} \\
 2-(24,24)-24 &  $\sqrt{2}\Lambda _{24} $ & $\sqrt{2}\Lambda _{24} $ & $P_{48n}$ &  \\
 2-(24,24)-24  & $\sqrt{2} O _{24} $ & $\sqrt{2} O _{24} $ & 
 $P_{48n}$, $P_{48p}$, $P_{48m}$ &   \\
\hline
\end{tabular} }
\end{table}

\begin{remark}
Table \ref{t1} lists the possible types of prime order automorphisms 
$\sigma \neq -1$. 
The type usually determines the genus of the $\Z $-lattice $L_K(\sigma )$
and the $\Z[\zeta _p]$-lattice $L_I(\sigma )$. 
If these genera are either classified in the 
literature (in particular the unimodular genera) 
or easily computed in Magma, we give the names or the
number of lattices of
minimum $\geq 6$ in these genera. 
Column ``example'' lists the known examples and the last column gives 
the two instances where the classification of the lattices 
is known to be complete. 
\end{remark}

\section{Automorphisms of type 
$5-(8,16)-8$.}

\begin{proposition}\label{unimodzeta}
There is a unique Hermitian unimodular $\Z[\zeta _5]$-lattice
$({\mathcal L},h)$ of dimension $\dim _{\Z[\zeta _5]} ({\mathcal L}) = 8$ 
such that the dual lattice 
$$\Lambda _{32} := ({\mathcal L},\frac{1}{5} \trace(h) ) ^{\#} $$ 
of the rescaled trace lattice of ${\mathcal L}$ has minimum $\geq 6$. 
\end{proposition}

\bew
A complete enumeration
of the genus 
using the Kneser neighboring method \cite{Kneser} which
is described for Hermitian lattices in \cite{Schiemann} shows that 
the genus 
of Hermitian unimodular $\Z[\zeta _5]$ lattices of rank 8 
contains 207 isometry classes of Hermitian unimodular lattices.
The completeness of the enumeration may also be verified using 
Shimura's mass formula, the mass is 
$$ 
\frac{46956347226527}{108864000000000} = \frac{46956347226527}{ 2^{15}3^55^97} 
\sim 
 0.43 .$$
Only for one of the 207 lattices the dual,
say $\Lambda _{32}$, of the trace lattice has
minimum $\geq 6$. 
\eb

The automorphism group of the $\Z $-lattice $\Lambda _{32}$ 
is a soluble group of order $2^83 \ 5^2 = 19200 $. 

\begin{theorem}\label{line9} 
There is a unique 
extremal even unimodular lattice of dimension $48$ that
admits an automorphism $\sigma $ of type $5-(8,16)-8$. 
The lattice is available under the name $P_{48m}$ in \cite{LAT}.
\end{theorem}

\bew 
Let $L$ be such a lattice and 
let $\sigma $ be an automorphism of Type $5-(8,16)-8$ of $L$.
Then $L_{K}(\sigma )$ is a 16-dimensional lattice of minimum 6 in the genus of the
5-modular lattices. It has been shown in \cite[Theorem 8.1]{BaVe} that this 
genus contains a unique lattice, say $\Lambda _{16}$, of minimum 6, which 
was denoted by $[2.\Alt_{10}]_{16}$ in \cite{NePl}.
\\
Since $s=z$ and by Remark \ref{unimod}, 
the dual $D = L_{I}(\sigma )^{\#} $ of the lattice $L_{I}(\sigma )$ is the trace lattice
of a Hermitian unimodular lattice over $\Z[\zeta _5]$ of dimension 8. 
So $L_I(\sigma ) \cong \Lambda _{32}$ by Proposition \ref{unimodzeta}. 
Therefore (up to isometry) $L$ contains a sublattice 
$M:= \Lambda _{32} \oplus \Lambda_{16}$ of index $5^8$. 

With extensive computations in {\sc Magma} \cite{MAGMA} 
%(rough guess of about 10 years of CPU time) 
we construct a set of lattices that
contains representatives of all $\Aut(M)$-orbits on 
$$ \{ X \mid M \leq X \leq M^{\#} , \min (X) =6 , X=X^{\#} \} .$$
To this aim we first compute orbit representatives of the 
1-dimensional subspaces of $\Lambda _{32}^{\#}/\Lambda _{32} $ under the
action of $\Aut(\Lambda _{32})$ to see that no proper integral overlattice of  
$\Lambda _{32}$ has minimum $\geq 6$. 
We then choose a suitable basis 
$\overline{b}:= (\overline{b}_1,\ldots , \overline{b}_8)$ of $\Lambda _{32}^{\#}/\Lambda _{32} $ 
such that $(b_i,b_i) = \frac{12}{5} = \min (\Lambda _{32} ^{\#} )$ and
$$(b_1,b_i) \in \frac{1}{5}+\Z , (b_2,b_j) \in \Z, 2\leq i\leq 8, 3\leq j\leq 8 .$$
More precisely the Gram matrix $(b_i,b_j)$ is
$$ F:= \left( \begin{array}{rrrrrrrr} 
12/5& 6/5& 6/5& 6/5& 6/5& 6/5& 6/5& 1/5 \\
 6/5&12/5&   1&   1&   1&   1&   1&   1 \\
 6/5&   1&12/5&   1& 1/5&   1& 1/5& 3/5 \\
 6/5&   1&   1&12/5& 4/5& 1/5& 2/5& 1/5 \\
 6/5&   1& 1/5& 4/5&12/5& 4/5& 4/5& 3/5 \\
 6/5&   1&   1& 1/5& 4/5&12/5& 2/5& 2/5 \\
 6/5&   1& 1/5& 2/5& 4/5& 2/5&12/5& 3/5 \\
 1/5&   1& 3/5& 1/5& 3/5& 2/5& 3/5&12/5 \end{array} \right) .$$
%Achtung, 3 hinter 6 eingefuegt, damit aus der Reihenfolge 4,5,6,3,7,8 
%die natuerliche Reihenfolge wird
Fixing the basis above  we hence know that 
$L = \langle M , (b_i, c_i) \mid 1\leq i\leq 8 \rangle $
where $(\overline{c}_1,\ldots , \overline{c}_8)$ is a basis of 
$\Lambda _{16}^{\#}/\Lambda _{16} $ so that 
that the Gram matrix of $c$ is congruent to $-F$ mod $\Z $.
Computing orbit representatives of the
action of $2.\Alt_{10}$ we 
find that all non-zero classes in $\Lambda_{16}^{\#}/\Lambda _{16}$
are represented by vectors of norm $\leq 4$.
In particular 
$\overline{c}_i = c_i + \Lambda_{16} $ has minimum $18/5$ for $i=1,\ldots , 8$. 
To narrow down the possibilities for $c_1$ and $c_2$ we 
use the fact that $\langle \overline{b}_1 , \overline{b}_2 \rangle $
contains six classes of minimum $12/5$. 
So we compute the orbit representatives of the
2-dimensional subspaces of $\Lambda_{16}^{\#}/\Lambda _{16} $
under  $2.\Alt _{10}$ to find that there are exactly 47 orbits of 
such 2-dimensional spaces that contain at least 6 classes of minimum 
$18/5$. These are our candidates for 
$\langle \overline{c}_1, \overline{c}_2 \rangle $. 
All of them contain exactly 6 such classes $(\pm \overline{x}, \pm \overline{y},\pm \overline{z})$, where the signs are chosen 
such that $(x,y) $ and $(x,z) $ are in $\frac{-1}{5} + \Z $. 

We now fix one of these 47 orbit representatives. 
As $\Aut(\Lambda _{32}) $ acts transitively on the minimal vectors of 
$\Lambda _{32}^{\#} $ we may assume that $c_1 = x$ is fixed. 
Then our candidates for $c_2 $ are $y$ and $z$. 
For ten of the spaces, the stabiliser of $x$ in the stabiliser of 
the 2-dimensional subspace interchanges $y$ and $z$, so we may assume that 
$c_2=y$ in these cases and the program splits naturally into 
84 subroutines each starting with a different $(c_1,c_2)$. 

For a fixed $(c_1,c_2)$ we compute the set $S(c_1,c_2)$
of all vectors $v$ of norm $18/5$
in $\Lambda _{16}^{\#}$ so that $(v,c_1) \in -\frac{1}{5} +\Z $ and 
$(v,c_2 ) \in \Z $. 
Then for $i=3,\ldots , 8$ the list of candidates for $c_i$ is 
$$C_i := \{ v\in S(c_1,c_2) \mid \min ( \langle M , (b_1, c_1), ( b_2 , c_2), 
(b_i, v) \rangle ) = 6 \} .$$
We get sets of different cardinalities varying between 
$1119 \leq |C_i| \leq 2469 $. 
Usually the smallest set is $C_3$ with a cardinality 
between 1119 and 1261. 
For all we successively construct 
$$L_k := \langle M, (b_1, c_1), (b_2 , c_2), 
(b_3 , c_3) , (b_4 , c_4),\ldots ,( b_k ,c_k ) \rangle $$
where the $c_i$ run through $C_i$ and we stop if 
$\min (L_k ) \leq 6$ or if $L_k$ is not integral. 
If we reach $k=8$ then we have found an extremal unimodular lattice. 
It takes between 1 and 2 hours to check all possibilities for a 
fixed $c_3$ (depending on the computer and on the deepness of the 
recursions). 
%Only in six cases, the computer checks the 
%case $k=7$, if one integral $L_7$ has minimum $\geq 6$ then 
%this branch continues to construct exactly one extremal $L_8$.
In total this algorithm constructs six  
 extremal even unimodular lattices $X$.

We check that all these six lattices 
 are indeed in the
same orbit under $\Aut (M)$. Let $P_{48m}$ be one representative 
of this orbit. 
Then $\Stab_{\Aut(M)} (P_{48m}) $ has order 1200 
and is identified with {\sc Magma} with group number 
$(1200,573)$. 
\eb

%As $\langle \sigma \rangle $ is the kernel of the action of $\Aut(M)$
%on $M^{\#}/M$ it is a normal subgroup of $\Aut (M)$. 
%So $$\Stab _{\Aut(M)} (P_{48m}) \leq N_{\Aut(P_{48m})} (\langle \sigma \rangle )
%=:N .$$
%On the other hand this  normaliser $N$ acts on $L_I(\sigma ) \oplus L_K(\sigma )=M$, \\
%so $N=\Stab _{\Aut(P_{48m})}(M)= \Stab _{\Aut(M)} (P_{48m}) $.
%
%For the natural epimorphisms $\pi _i : \Aut(M) \to \Aut(\Lambda _i)$ 
%($i=32,16$) we get that $\pi_{32} (N) \cong N$ and $\pi_{16} (N) \cong N/\langle \sigma \rangle $. 

\begin{remark}
The normalizer $N$ of $\langle \sigma \rangle$ in $\Aut(P_{48m})$ 
has order $1200$ and is isomorphic to
group number $(1200,573)$ in the {\sc Magma} small groups library.
The Sylow 2-subgroups of $\Aut(P_{48m})$ are isomorphic to the 
ones in $N$, $D_8 {\sc Y} C_4 \in \Syl_2(P_{48m})$
and the Sylow 5-subgroups of $\Aut(P_{48m})$ are 
isomorphic to $C_5\times C_5 \in \Syl _5(N)$.
\end{remark}

\bew
Let $L := P_{48m}$, 
$G:=\Aut(L)$ and $N:= N_G(\langle \sigma \rangle )$.
Then $N$ acts on $M:=L_{K}(\sigma ) \operp L_I(\sigma ) $, 
so $N= \Stab _{\Aut(M)} (P_{48m}) $ is as computed in the proof of
Theorem \ref{line9}.
Let $S$ be a Sylow 2-subgroup of $N$. 
With {\sc Magma} we compute that
$M:=\sum _{s\in S} (s-1) L$ is a sublattice of $L$ of index $2^6$,
minimum 6 and kissing number 821760.
{\sc Magma} computes that the automorphism group of $M$ is equal to $S$.
As $N_G(S)$ acts on $M$ we conclude that $N_G(S) = S$.
As proper subgroups of $p$-groups are strictly contained 
in their normalizers, this implies that $S$ is a
Sylow 2-subgroup of $G$.
Let $S$ be a Sylow 5-subgroup of $N$.
As $\langle \sigma \rangle $ contains all elements of 
order 5 with a 16-dimensional fixed space, we get that 
$N_G(S) \subseteq N_G(\langle \sigma \rangle )  = N$. 
In particular $S\in \Syl_5(G) $.
\eb

\end{document}